\documentclass[12pt]{amsart}

\usepackage{url}

\usepackage{graphicx}

\theoremstyle{plain}
\newtheorem{lem}{Lemma}
\newtheorem{thm}{Theorem}
\newtheorem*{thm*}{Theorem}
\newtheorem{cor}{Corollary}

\newcommand{\boundary}{\partial}
\newcommand{\coboundary}{\delta}
\newcommand{\equivalent}{\sim}
\newcommand{\Id}{\mathbf1}

\newcommand{\wm}{\varepsilon}   
\newcommand{\p}{\hspace*{-3pt}+}
\newcommand{\n}{\hspace*{-3pt}-}
\newcommand{\Z}{\mathbf Z}

\title{Cohomology and Immersed Curves}
\author{Mario O. Bourgoin}
\date{\today}
\address{Department of Mathematics \\
Brandeis University \\
415 South Street \\
MS 050 \\
Waltham, MA 02454}
\email{mob@brandeis.edu}

\begin{document}

\begin{abstract}
  We introduce a new cohomology-theoretic method for classifying
  generic immersed curves in closed compact surfaces by using Gauss
  codes.  This subsumes a result of J.S. Carter on classifying
  immersed curves in oriented compact surfaces, and provides a
  criterion for when an immersion is $2$-colorable.  We note an
  application to twisted virtual link theory.
\end{abstract}

\maketitle

\section{Introduction}

\begin{figure}[htbp]
  \centering
  \includegraphics[height=1.5in]{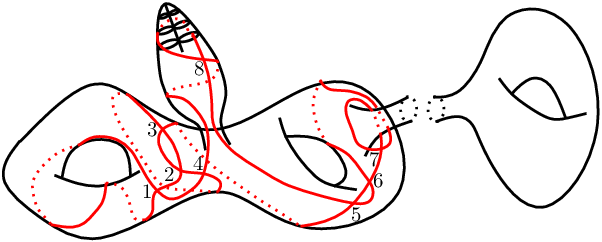}
  \caption{An immersed curve with (de)stabilization.}
  \label{fig:curve1}
\end{figure}

We associate with each generic immersed curve in a closed compact
surface an isomorphism class of \emph{intersigned Gauss codes}, which
are multi-component double occurrence sequences of symbols drawn from
a fixed alphabet whose symbols alternate with signs, with the ordered
signs called a \emph{sign sequence}.  Figure~\ref{fig:curve1} shows an
example of such a curve in a non-orientable surface,
and a representative of its intersigned Gauss code is
\begin{align*}
  1+2-3-4-2-3+5-6+7-7+6-5+8+8-4+1-
\end{align*}
We have the following.
\begin{thm}\label{thm:classify_curves}
  Stable geotopy classes of generic immersed curves correspond to
  isomorphism classes of intersigned Gauss codes.
\end{thm}

When the generic immersion of curves in a surface is cellular, it
determines an embedded graph which is the $1$-skeleton of a cellular
decomposition of the surface, and then the sign sequence of the
intersigned Gauss code is interpreted as a $1$-cochain (i.e. an
assignment of $\pm1$ to each of the edges) on this $1$-skeleton.  By
pulling back the intersection points of the curve to the circles in
the domain of the immersion, we can see that the embedded graph is a
quotient of the union of cycle graphs, and the sign sequence is a
cocycle on this graph.  There is an obvious correspondence between the
cycle graphs with a cocycle and the intersigned Gauss code.

\begin{figure}[htbp]
  \centering
  \raisebox{-.725in}{\includegraphics[height=1.55in]{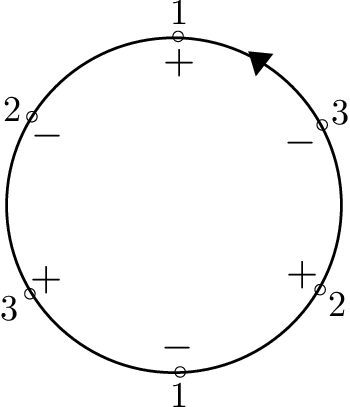}}~
  $\buildrel K\over\Longrightarrow$~
  \raisebox{-.625in}{\includegraphics[height=1.25in]{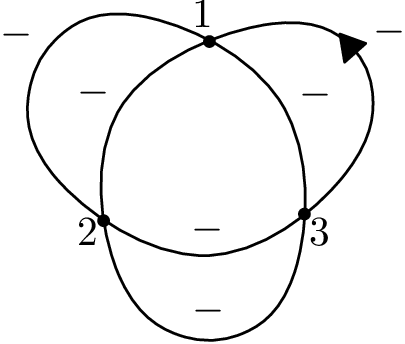}}
  \caption{A curve represented by $1\n2\n3\n1\n2\n3\n$.}
  \label{fig:cycle}
\end{figure}

\begin{figure}[htbp]
  \centering
  \includegraphics[height=3in]{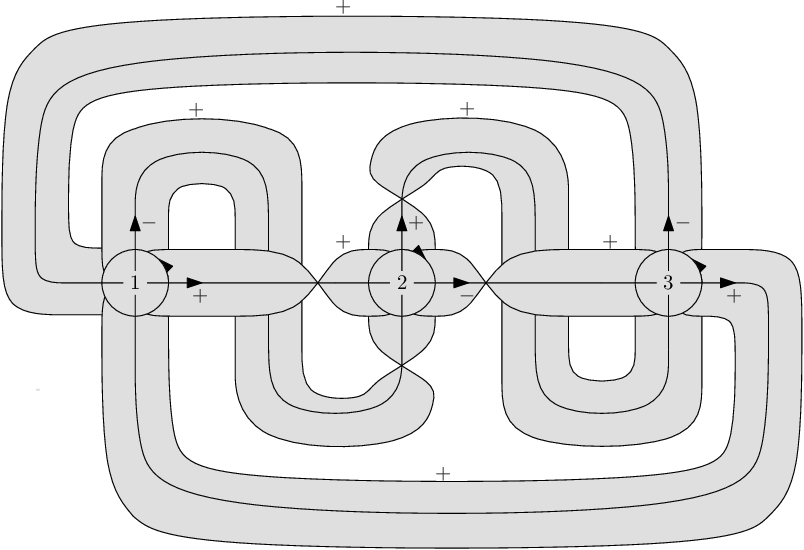}
  \caption{The immersed curve of Figure~\protect\ref{fig:cycle} in a ribbon graph.}
  \label{fig:sphere}
\end{figure}

A \emph{rotation cochain} is a $0$-cochain on a cycle graph that
assigns opposite signs to the two pullbacks of the intersection
points.  This situation is illustrated in Figure~\ref{fig:cycle} where
a cycle graph with a rotation cochain is shown on the left, and the
associated immersion in the $2$-sphere is shown with the associated sign
sequence on the right.
It is clear that the rotation cochain also defines an assignment of
opposite signs to the symbols of the associated intersigned Gauss
code.  Then we have the following:
\begin{thm}\label{thm:orientable_surfaces}
  Let $C$ be a curve cellularly immersed in a closed surface $\Sigma$,
  $\Pi$ be an intersigned Gauss code for $C$, and $r$ a rotation
  cochain on $\Pi$.  The sign sequence of $\Pi$ added mod~$2$ to the
  coboundary of $r$ is a cocycle representing the first
  Stiefel-Whitney class $w_1$ of $\Sigma$.
\end{thm}
This has the immediate corollary:
\begin{cor}
  The closed surface $\Sigma$ is orientable if and only if the sign
  sequence of $\Pi$ plus the coboundary of $r$ is a mod~$2$ coboundary.
\end{cor}
Figure~\ref{fig:sphere} shows a ribbon graph containing the oriented
immersed curve of Figure~\ref{fig:cycle}.
As the ribbon graph has five boundary components, the closed surface
obtained from it by attaching a disk to each component has Euler
characteristic $2$ showing it is the $2$-sphere.  The curve's
orientation is shown by arrowheads on the signed curve segments
through the crossings; a segment's sign is shown next to its
arrowhead.  The figure shows the orientations of disk-shaped
neighborhoods of the curve's three crossings with respect to the plane
of the paper by the arrowheads on the disk boundaries.  This
orientation is determined by the direction of rotation of the
crossing's $+$ segment about the crossing to its $-$ segment without
passing through the $-$ segment so that their arrowheads overlap.  The
assignment of the arrowhead signs is in turn is set to reproduce the
sign sequence of the rotation cochain on the circle graph on the left
side of Figure~\ref{fig:cycle}.  The sign shown along a ribbon
connecting two disks determines whether the ribbon twists to match
($+$) or not ($-$) the orientations the end disks as either disk's
orientation is propagated through the ribbon to the other disk.  The
cochain that assigns those signs to ribbons is a cocycle representing
of $w_1$.  It is calculated by ``adding'' in the usual way the signs
of the curve segments the ribbons connect to the respective intersigns
between the crossings; the calculation is shown in Table~\ref{tab:w1}.
\begin{table}[htbp]
  \centering
  \begin{tabular}{r|cccccccccccc}
    Intersigned Gauss Code & 1 & $-$ & 2 & $-$ & 3 & $-$ & 1 & $-$ & 2
    & $-$ & 3 & $-$ \\
    Ribbon Start && $+$ && $-$ && $+$ && $-$ && $+$ && $-$ \\
    Ribbon End && $-$ && $+$ && $-$ && $+$ && $-$ && $+$ \\
    \hline
    Cocycle && $+$ && $+$ && $+$ && $+$ && $+$ && $+$ \\
  \end{tabular}
  \caption{Calculating the $w_1$-cocycle for the curve of Figure~\protect\ref{fig:cycle}}
  \label{tab:w1}
\end{table}

As to the $2$-colorability of a cellular immersion, we have the
following:
\begin{thm}\label{thm:two-colorable_immersion}
  Let $C$ be a curve cellularly immersed in a closed surface $\Sigma$
  and $\Pi$ be an intersigned Gauss code for $C$.  The mod~$2$
  complement of the sign sequence of $\Pi$ is a cocycle representing
  the mod~$2$ Poincar\'e dual of the mod~$2$ homology class of $C$.
\end{thm}
Then the theorem has the immediate corollary:
\begin{cor}
  The immersion $C$ is $2$-colorable if and only if the mod~$2$
  complement of the sign sequence of $\Pi$ is a mod~$2$ coboundary.
\end{cor}

\section{Background}

\subsection{Gauss Codes}

Gauss introduced the concept of a Gauss code to describe a normal
immersion of a circle in the plane~\cite[pp. 271--286]{MR82e:01122h}.
Such a code is obtained from an immersion by labeling its crossings
and reading them off starting at an arbitrary non-crossing point on
the curve.

A \emph{Gauss code} is a finite collection of (non-empty) finite
sequences of symbols drawn from a set such that every symbol of the
set occurs exactly twice in the sequences of the collection.  Two
Gauss codes over the same symbols are equivalent under permutation of
the set of symbols, and rotation and reflection of the sequences.
Each finite sequence is a \emph{component} of the Gauss code.

\begin{sloppypar}
  Gauss codes with one component are sometimes called Gauss
  words~\cite{MR91d:57002}, double occurrence words or
  sequences~\cite{MR58:10883}, and cross codes.  And when
  multi-component Gauss codes have been used to classify a normal
  immersion of a disjoint union of circles in an oriented closed
  compact surface, they have been called Gauss
  paragraphs~\cite{MR91d:57002}.
\end{sloppypar}



A Gauss code does not uniquely specify a collection of planar curves
as there are two ways of realizing each crossing in the immersion
surface.  If the curves are oriented, then at every crossing, one
strand is traversed left to right by the other strand, while the other
strand is traversed right to left by the first strand.  If the
orientation of the components is fixed and opposite signs are
associated with each symbol of the code to indicate the direction of
traversal, the Gauss code is called a \emph{signed Gauss code}, an
\emph{intersection sequence}~\cite{MR0252254}, or a
\emph{lacet}~\cite{MR2002e:05048}.  Signed Gauss codes were used by
J.S.~Carter~\cite{MR91d:57002} to classify stable geotopy classes of
generic immersions of curves in closed oriented surfaces of higher
genus.

As Gauss pointed out, some Gauss codes do not have an associated
planar diagram.  For example, the Gauss code $1212$ is not planar.
Gauss identified a necessary condition for a Gauss code to be that of
a planar diagram, namely that there be an even number of symbols
between the two occurrences of any single symbol in the code.  This
condition is sometimes called \emph{evenly intersticed} or
\emph{evenly interlaced}.  That this condition is not sufficient to
characterize planar Gauss codes is shown by the Gauss code
$1234534125$.  Gauss conjectured that the interlacement structure of a
Gauss code would determine whether it has a planar diagram.  This
conjecture was first proved by Rosenstiehl~\cite{MR58:10883}.

\subsection{Lacets}

One direction of development in the study of Gauss codes came from
Lins, Richter, and Shank~\cite{MR88i:05071} as an extension of the
search for combinatorial conditions for the planarity of Gauss codes.
In that paper, the authors developed an algebraic approach to
determine the possible surfaces in which a given Gauss code is
associated with a $2$-colorable immersed curve based on a partition of
the symbols of the Gauss code.  These $2$-colorable immersions of
circles have been called \emph{lacets}~\cite{MR2002e:05048}.
Cohomology-theoretic obstructions have been developed to determine
conditions under which a Gauss code can be presented as a lacet in a
particular surface.  These depend upon the structure of the
\emph{interlacement graph} of a Gauss code.  This graph has for
vertices the symbols of the Gauss code, and has an edge between two
vertices whose symbols are interlaced in the Gauss code.  Two symbols
$1$ and $2$ in a Gauss code are interlaced if they appear in a pattern
of the form $\cdots1\cdots2\cdots1\cdots2\cdots$.  Typical results
include the following theorem from~\cite{MR88i:05071}.
\begin{thm*}[Lins et al.~(1987)]
  The $2$-colorable immersions of a Gauss code are in an orientable
  surface if and only if every vertex of its interlacement graph is
  even-valent.
\end{thm*}
So the Gauss code ``$1234534125$'' has $2$-colorable immersions in
orientable surfaces while the Gauss code ``$1212$'' has $2$-colorable
immersions in non-orientable surfaces.

\section{Definitions}

A \emph{generic immersion} of curves is the immersion of a disjoint
union of circles in a surface such that the only intersections are
transverse double points.  Two generic immersion of curves are
\emph{geotopic} if there is a homeomorphism between the immersion
surfaces that takes the image of one immersion onto the image of the
other immersion.  The immersions are \emph{stably geotopic} if there
is a collection of $1$-handles and crosscaps that can be added to (or
removed from) either surface in the complement of the immersions such
that the resulting immersions are geotopic.  An immersion is
\emph{cellular} if its complement is homeomorphic to a collection of
open disks.  It is \emph{$2$-colorable} if the components of its
complement can be assigned one of two colors such that the arcs
between intersections always separate one color from the other.  Thus
it is $2$-colorable if and only if it represents a null-homologous
cycle mod~2.

A \emph{generic immersion} of oriented curves is a generic immersion
of curves where the circles are oriented, and then two generic
immersions of oriented curves are \emph{geotopic} when the
homeomorphism preserves the orientations.

\emph{Signed Gauss codes} are presentations of permutations $P$ of the
signed symbols $\{1^\pm, \ldots, n^\pm\}$ over an alphabet $\{1,
\ldots, n\}$.  The \emph{components} of a signed Gauss code are the
orbits of $P$.  Two signed Gauss codes are \emph{isomorphic} if one
can be obtained from the other by a cyclic permutation of any
component, by permuting the order of the components, by permuting the
alphabet, by changing all signs to their opposite, or by reversing any
component while complementing the signs on all symbols that occur on
two components where only one component is reversed.  For example, the
three component signed Gauss code
\begin{align*}
  &1^-2^+3^-1^+/2^-/3^+ &&\text{is isomorphic to} &&1^+3^+2^-1^-/3^-/2^+.
\end{align*}

\emph{Oriented intersigned Gauss codes} are multi-component double
occurrence sequences over an alphabet $\{1, \ldots, n\}$ whose symbols
alternate with signs, where the sequence of signs is called the
\emph{sign sequence}.  Two oriented intersigned Gauss codes are
\emph{isomorphic} if one can be obtained from the other by a cyclic
permutation of any component, by permuting the order of the
components, or by permuting the alphabet.  For example, the two
component intersigned Gauss code
  \begin{align*}
  &1-2+1-3-/2+3- &&\text{is isomorphic to} &&2+3-/2+1-3-1-.
\end{align*}

\emph{Intersigned Gauss codes} are oriented intersigned Gauss codes
that admit the additional isomorphism of reversing any component while
at the same time complementing the signs between all pairs of symbols
where both occurrences of one symbol are on the reversed component and
only one occurrence of the other symbol is on that component.  For
example, the two component intersigned Gauss code
  \begin{align*}
  &1-2-3+1-/2+3- &&\text{is isomorphic to} &&1-3-2+1-/2+3-.
\end{align*}
In this transformation, the first component is reversed.  At the same
time, since on that component $1$ occurs twice and $2$ occurs once,
the sign between $1$ and $2$ is changed from $-$ to $+$.  And for
similar reasons, the sign between $1$ and $3$ is changed from $+$ to
$-$.


Unless otherwise said, graphs may have multiple edges and loops.  An
\emph{Euler partition} of a graph is a collection of closed paths that
together go over every edge exactly once.  A graph admits an Euler
partition if and only if it is even-valent.  An embedding of a graph
in a surface that may or may not be orientable is an embedding of the
corresponding $1$-complex, and induces a cyclic order of the edges at
every vertex.  An embedding of a graph that has an Euler partition is
\emph{straight-through at a vertex with respect to that partition} if
for all edge pairs determined at that vertex by the closed paths, the
two edges have the same number of other edges between them in either
direction of the cyclic order of the edges at that vertex induced by
the embedding of the graph.  An embedding is called just
\emph{straight-through} if it is straight-through at every vertex.  An
embedding that is not straight-through at a vertex is \emph{bent} at
that vertex.

Henceforth, we restrict ourselves to $4$-valent graphs to simplify the
exposition.

An \emph{embedding scheme} for a graph consists of a pair $(\pi,
\lambda)$ of a \emph{rotation system} $\pi$ which is a set of cyclic
permutations $\pi_v$ of the edges incident with a vertex $v$, and a
\emph{signature} $\lambda$ from the set of edges $e$ to
$\pm1$~\cite{MR2002e:05050}.  Each embedding scheme for a graph
encodes a \emph{cellular} immersion for the graph, where the
complement of the image of the graph is homeomorphic to a collection
of open disks.

The map between embedding schemes and cellular immersions is as
follows.  Disjoint neighborhoods of each vertex in the graph are
chosen.  The neighborhood of each vertex is embedded in an oriented
disk with boundary such that the counterclockwise order of the
vertex's incident edges on the boundary is given by its permutation.
For each disk, disjoint half-disk neighborhoods of the intersection of
an edge with the boundary of the disk are chosen, and obtain their
orientation from that of the disk.  Each pair of neighborhoods
intersecting an edge are identified by orientation-preserving
(resp.~orientation-reversing) homeomorphisms if the edge's signature
is $+$ (resp.~$-$).  This extends the embedding of the vertices in
disks to the graph in a compact surface.  The resulting surface is
closed by attaching disks to its boundary's components.

Conversely, given an embedded graph in a closed compact surface,
choose a collection of disjoint oriented closed neighborhoods
homeomorphic to disks about the vertices to determine a rotation
system.  The signature is determined by propagating the orientations
of the disks along tubular neighborhoods of the edges, and is $+$ when
the orientations agree and $-$ otherwise.  An edge is \emph{twisted}
with respect to an embedding scheme if its signature is $-$.  Every
cellular embedding of a graph in a surface is uniquely determined, up
to homeomorphism, by an embedding scheme.  Two embedding schemes for a
graph are \emph{equivalent} if one can be obtained from the other by
reversing some of the cyclic permutations along with inverting the
signatures on their incident edges.  Equivalence classes of embedding
schemes for two graphs are \emph{isomorphic} if there is an
isomorphism between the graphs that respects the rotation system and
signature of representative embedding schemes.

\section{The Classification of Immersed Curves}

We begin by classifying cellular immersions of oriented curves, then
we classify cellular immersions of curves, and finally we use
surgeries to obtain a classification for all generic immersed curves.

\begin{lem}\label{lem:ic_es}
  Geotopy classes of generic cellularly immersed curves correspond to
  isomorphism classes of straight-through embedding schemes.
\end{lem}
\begin{proof}
  Generic immersions of unoriented curves can be associated with
  four-regular embedded graphs in an obvious way, so all generic
  cellularly immersed curves can be determined from the equivalence
  classes of embedding schemes of their graphs.  A graph that comes
  from a generic immersion of curves has associated with it an Euler
  partition.  Since our immersions are transverse, we need only
  consider embedding schemes that are straight-through with respect
  to the Euler partition.
\end{proof}

\begin{figure}
  \centering
  \includegraphics[width=\textwidth]{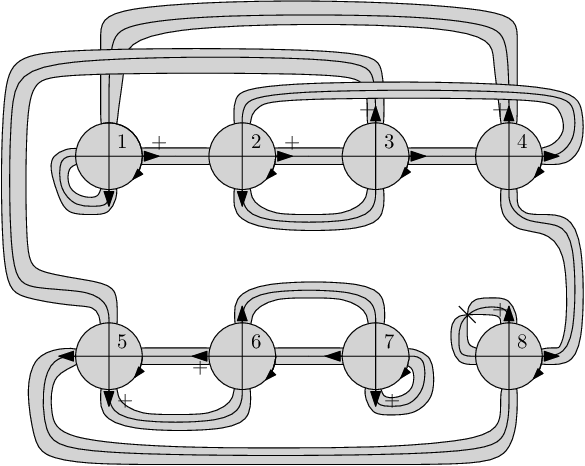}
  \caption{The immersed curve of Figure~\ref{fig:curve1} in a ribbon graph}
  \label{fig:curve2}
\end{figure}

Since the embedding schemes are straight-through, we have two choices
of permutation for each intersection and two choices of sign for each
edge.  If we fix the choices of rotations for all intersections of
generic immersions of curves that have isomorphic graphs, we can
determine a cellular generic immersion of curves by the choice of edge
signs.  The problem is then to find a presentation of the class of
embedding schemes for generic immersions of curves that is independent
of the choice of rotation system.
\begin{lem}\label{lem:geo_unique}
  Geotopy classes of generic cellularly immersed oriented curves
  correspond to isomorphism classes of oriented intersigned Gauss
  codes.
\end{lem}
Figure~\ref{fig:curve2} is the immersed curve of
Figure~\ref{fig:curve1} as a ribbon graph.  For convenience, its
crossings are in disks that are embedded in the plane of the paper and
all given the same clockwise orientation that turns the $+$ arrow head
into the $-$ arrow head.  Then, the ribbon edges all connect
neighborhoods with the same orientation.  In Figure~\ref{fig:curve1},
after the curve passes through crossing $8$, it goes through a
crosscap, and then passes through crossing $8$ again.  Then, the curve
in the ribbon graph must be in a twisted ribbon edge to change the
local orientation along the curve.  The ribbon graph has $7$ boundary
components so if we closed the surface using disks, the resulting
non-orientable surface would have Euler characteristic $-1$.
\begin{proof}[Proof of lemma~\ref{lem:geo_unique}]
  Given an oriented intersigned Gauss code, define a union of cycle
  graphs with two distinct vertices for each symbol in the alphabet,
  and an edge for each pair of cyclically adjacent symbols in some
  component, where the edge goes between the distinct vertices
  associated with the symbols.  Define a four-regular graph as the
  quotient of the cycle graphs that identifies each pair of distinct
  vertices that correspond to a symbol.  Since the edges of the cycle
  graphs are in $1$-to-$1$ correspondence with those of the
  four-regular graph, by abuse of notation we will identify the edges
  of the cycle graphs and those of the four-regular graph.

  We define an embedding scheme for the four-regular graph by first
  choosing a rotation system that is straight-through at every vertex
  with respect to the Euler partition of the graph given by the
  intersigned Gauss code.  This defines a rotation cochain on the
  cycle graphs that assigns a $+$ to the pullback of a vertex whose
  outgoing edge is sent by the rotation system to the vertex's other
  outgoing edge, and a $-$ to the other pullback of the vertex.  Since
  the sign sequence of the intersigned Gauss code is a $1$-cocycle of
  the cycle graphs, we add it mod~$2$ to the coboundary of the
  rotation cochain and define it to be the embedding scheme's
  signature.  Since the mod~$2$ sum of two rotation cochains is the
  pullback of a $0$-cochain on the embedded graph, any another choice
  of rotation cochain would produce an equivalent embedding scheme.
  Any other isomorphic oriented intersigned Gauss code will have an
  isomorphic graph, and since the above map does not depend upon the
  presentation of the Gauss code, it will have an isomorphic embedding
  scheme.

  Conversely, given a generic cellular immersion of oriented curves in
  a surface, label the $n$ intersections with symbols from an alphabet
  $\{1, \ldots, n\}$.  By Lemma~\ref{lem:ic_es} we can obtain a
  representative of an equivalence classes of embedding schemes.  The
  choice of rotation system along with the orientation of the curves
  determines a $+$ sign for the occurrence of the symbol whose
  outgoing edge is sent to the other occurrence's outgoing edge, and a
  $-$ sign for the other occurrence.  This is a rotation cochain.  We
  can then read off an intersigned Gauss code whose sign sequence is
  the mod~$2$ sum of the embedding scheme's signature and the
  coboundary of the rotation cochain.

  This is illustrated in Figure~\ref{fig:signature} where we obtain
  the sign between symbols $1$ and $2$ to be a $-$.  Each symbol
  identifies a vertex of the graph that is embedded in an oriented
  neighborhood, itself contained in the plane of the sheet on which
  the figure is drawn; assuming the sheet has the standard
  counter-clockwise orientation, vertex $1$'s neighborhood is oriented
  clockwise while vertex $2$'s neighborhood is oriented
  counter-clockwise.  The edge joining the vertices passes through an
  attaching band that matches the vertex neighborhoods' orientations;
  because these neighborhoods have opposite orientations, the
  attaching bad twists to match them, and so cannot lie in the plane
  of the figure's sheet.  The sign between $1$ and $2$ is then the sum
  of the sign of the edge's segment in the oriented neighborhood of
  vertex $1$, a $-$, the sign of the edge's segment in the oriented
  neighborhood of vertex $2$, a $+$, and the sign of the edge's
  segment connecting them that is embedded in an orientation-matching
  neighborhood, a $+$.\footnote{The signs $+$ and $-$ behave as,
    respectively, the $0$ and $1$ of the ring $\Z/2\Z$.}
  \begin{figure}[htbp]
    \centering
    \includegraphics[width=\linewidth]{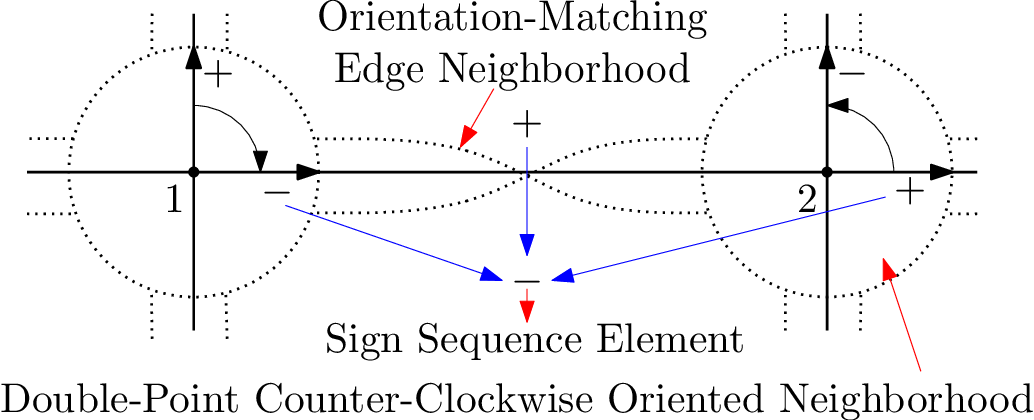}
    \caption{Calculating the sign between two symbols.}
    \label{fig:signature}
  \end{figure}

  Changing the choices of assignment of symbols to intersections,
  representative embedding scheme, starting point on the curves, and
  order of curves in reading the intersigned Gauss code will produce
  an isomorphic oriented intersigned Gauss code.
\end{proof}

We characterize the differences between the intersigned Gauss codes
for generic cellular immersions of oriented curves that differ only in
the orientations of their curves.
\begin{lem}
  Two generic cellular immersions of oriented curves that differ only
  in the orientation of one component have oriented intersigned Gauss
  codes that differ in the reversal of the corresponding component
  while changing the signs between all pairs of symbols where one
  symbol occurs on a reversed component and the other does not.
\end{lem}
\begin{proof}
  Reversing only the direction on a component while keeping fixed
  the other choices will assign the opposite straight-through
  rotation to the intercomponent intersections on that component,
  and this will change the signs between those intersections and
  other intersections.
\end{proof}

It immediately follows that geotopy classes of generic cellularly
immersed curves correspond to isomorphism classes of intersigned Gauss
codes.

We consider general generic immersions of curves.
\begin{lem}\label{lem:unique}
  Stable geotopy classes of generic immersions of curves have a unique
  cellular immersion.
\end{lem}
\begin{proof}
  The complement of the image of a generic immersion of oriented
  curves in a surface is homeomorphic to a set of open surfaces, and
  the individual closure of each component has a collection of
  circles for a boundary.  Then, a collection of the same number of
  open disks may be substituted for the surface components giving
  the unique cellular immersion.
\end{proof}

\begin{proof}[Proof of Theorem~\ref{thm:classify_curves}.]
  By Lemma~\ref{lem:unique}, any generic immersion of curves has a
  unique cellular immersion from which we may obtain an isomorphism
  class of intersigned Gauss codes.  Since the stabilizations occur on
  the complement of the immersed curve, they do not change the
  isomorphism class of intersigned Gauss codes.
\end{proof}

\section{Orientable Surfaces}

One lemma is needed:
\begin{lem}
  The signature of an embedding scheme is a representative of $w_1$ of
  the associated surface.
\end{lem}
\begin{proof}
  Fix an embedding scheme for a given surface.  The signature of the
  embedding scheme is a $1$-cochain in the cellular cohomology of the
  surface that, along with the rotation scheme, determines the face
  cycles.  A $1$-cochain is a cocycle of the surface if and only if it
  has an even number of negative edges on every face cycle.  Since the
  face cycles are the boundaries of disks, they must be
  orientation-preserving, and so have an even number of twisted edges
  along their path.  Obtain another embedding scheme from the surface
  by choosing the same orientation for every disk containing a vertex,
  so that the twisting of the edges determines the signature of the
  scheme.  For this scheme, there is an even number of negative edges
  along each face cycle determined by the signature.  The two
  embedding schemes differ by the choice of the permutation that
  represents the order of the edges about the vertex, and where they
  differ defines a $0$-cochain of the surface.  Then, the signatures
  of the two embedding schemes differ by the coboundary of this
  $0$-cochain, so the original signature also had an even number of
  negative edges along each face cycle.

  By cellularity, any cycle in the surface is homologous to a cycle
  on the immersed curve, and a cycle in the immersed curve is
  orientation-preserving if and only if it has an even number of
  twisted edges along its path.  But this is true if and only if the
  path has an even number of negative edges determined by the
  signature.
\end{proof}

The proof of the characterization of the orientability of the
immersion surface is now straightforward.
\begin{proof}[Proof of Theorem~\ref{thm:orientable_surfaces}]
  Since the signature of an embedding scheme is a representative of
  $w_1$ of the surface, then the surface is orientable if and only if
  the signature is a coboundary.  This is true if and only if the
  mod~$2$ sum of the sign sequence of the intersigned Gauss code with the
  coboundary of a rotation cochain is a coboundary of a pullback
  cochain, or in other words if the sign sequence of the intersigned Gauss
  code is a coboundary of the mod~$2$ sum of a rotation cochain and a
  pullback cochain.  Since this sum is itself a rotation cochain, we
  are done.
\end{proof}

We now provide a map from the signed Gauss codes used by Carter and
intersigned Gauss codes.  Clearly, the signs of signed Gauss codes
define a rotation cochain on the associated cycle graph.  Then, the
sign sequence of the corresponding intersigned Gauss code is the
coboundary of that rotation cochain.  By
Theorem~\ref{thm:orientable_surfaces}, the associated surface is
orientable, and by construction, the associated cellular immersion is
geotopic.

\section{$2$-Colorable Immersions}

In this section, unless otherwise indicated, all graphs are
four-valent embedded graphs, except that the lifts of such graphs are
lifts to cycle graphs whose components map to straight-through closed
paths in the graphs.

\begin{proof}[Proof of Theorem~\ref{thm:two-colorable_immersion}]
  To show that the mod~$2$ complement of the signature is a cocycle
  representing the mod~$2$ Poincar\'e dual of the mod~$2$ homology class
  of $C$, we use the fact that if a cohomology class represented by a
  cocycle $Z$ is the dual of the homology class represented by a cycle
  $z'$, then
  \begin{align*}
    Z(z) &= z'\cdot z,~\text{for all cycles}~z.
  \end{align*}
  By cellularity, any cycle in the surface is homologous to a cycle on
  the immersed curve.  We associate a cycle on the immersed curve with
  a homologous cycle chosen such that:
  \begin{itemize}
  \item it is an $\epsilon>0$ push-off from the curve in the direction
    of a vector translated along the curve,
  \item intersects with the immersed curve once at every
    straight-through crossing,
  \item intersects with the immersed curve an even number of times at
    a bent crossing, and
  \item intersects the immersed curve one final time if the cycle is
    orientation-reversing.
  \end{itemize}
  This is illustrated in Figure~ref{fig:Candc}.
  \begin{figure}[htbp]
    \centering
    \begin{tabular}{ccc}
      Bend & Bend & Straight-Through \\
      \includegraphics[height=1.5in]{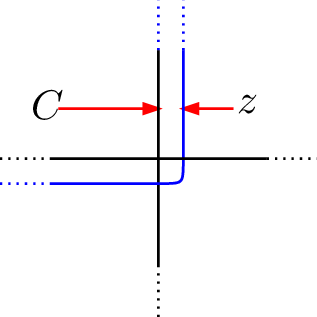} &
      \includegraphics[height=1.5in]{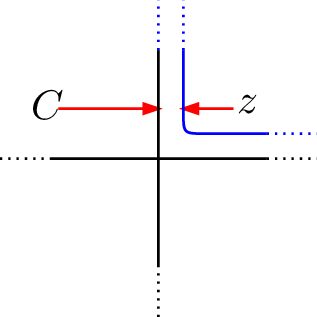} &
      \includegraphics[height=1.5in]{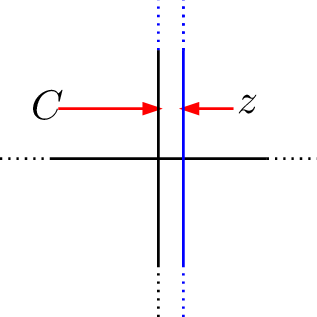}
    \end{tabular}
    \caption{Choosing a homologous cycle.}
    \label{fig:Candc}
  \end{figure}
  Then for $C$ being the immersed curve by a map $L$ and $z$ a cycle
  that meets the above conditions, we have that:
  \begin{align*}
    C \cdot z &= \text{\#straight-throughs}(z) + w_1(z) \pmod 2.
  \end{align*}
  Let $r$ be a rotation cochain, $S(\Pi)$ be the $1$-cochain
  determined by the sign sequence of $\Pi$, and $\Id$ be the
  constant $1$-cochain that maps each edge to $1$.  We have:
  \begin{align*}
    S(\Pi) + \Id &= (S(\Pi) + \coboundary r) + (\Id +
    \coboundary r) \pmod 2,\\
    \intertext{and we know from the proof of
      Theorem~\ref{thm:orientable_surfaces} that $[S(\Pi) +
      \coboundary r] = w_1$.  It is clear that:}
    \Id(z) &= \#\text{arcs}(z) \pmod 2. \\
    \intertext{By identifying the arcs of the embedded curve with
      those of the pullback by $L$ to the domain circles, we have:}
    \coboundary r(L^{-1}(z)) &= r(\boundary L^{-1}(z)) =
    \#\text{bends}(z) \pmod 2 \\
    \intertext{because $r$ has opposite values on the two lifts of a
      double point.  Then:} \Id + \coboundary r &=
    \#\text{straight-throughs} \pmod 2.
  \end{align*}
  Since the faces of an immersion are orientation-preserving paths
  that bend at every crossing, then $S(\Pi) + \Id$ is a cocycle, and
  it represents the mod~$2$ Poincar\'e dual of $C$.
\end{proof}

\section{Related Work}

As shown in~\cite{bourgoin08:_twist_link_theor}, when the immersion
surface of a curve is thickened by its orientation $I$-bundle, it
becomes possible to resolve the intersection points of the curve to
obtain a link in an orientable three-manifold, and then the immersed
curve becomes a diagram of this link.  We have classified these
diagrams using \emph{intersigned link codes}, which are complemented
intersigned Gauss codes with the addition of a writhe at every
crossing.  A \emph{Reidemeister move} on an intersigned link code is
one of the abstractly-defined moves:
\begin{center}
  \begin{align*}
    &\text{R-1}: &&ab~\equivalent~a 1^\wm\p1^\wm b \\
    &\text{R-2}: &&a 1^\wm\p2^{-\wm} b 1^\wm\p2^{-\wm}
    c~\equivalent~abc~\equivalent~a 1^\wm\p2^{-\wm} b 2^{-\wm}\p1^\wm c
    \\
    &\text{R-3}: &&a 1^\wm\p2^\wm b 3^\wm\n2^\wm c 1^\wm\n3^\wm
    d~\equivalent~a 1^\wm\p2^\wm b 1^\wm\n3^\wm c 3^\wm\n2^\wm d \\
    &&&a 1^{-\wm}\p2^\wm b 3^\wm\n2^\wm c 1^{-\wm}\n3^\wm d~\equivalent~a
    1^\wm\p2^{-\wm} b 1^\wm\n3^\wm c 3^\wm\n2^{-\wm} d
  \end{align*}
\end{center}
In a move, intersigned link codes are presented left-right, the
crossing numbers are assigned to reflect the order in which crossings
are encountered, the writhe mark is given as an exponent $\wm = +$ or
$-$ such that $-\wm =$, respectively, $-$ or $+$, and the lowercase
letters $a, b, c, \ldots$ represent segments of the code in between
which the fragments are embedded.  Two intersigned link codes are
\emph{Reidemeister equivalent} if there is a sequence of Reidemeister
moves or intersigned Gauss code isomorphisms that takes one code to
the other.  The resulting theory of twisted virtual links is a proper
extension of both Lou~Kauffman's virtual link
theory~\cite{MR2000i:57011} and Yu.~V.~Drobotukhina's projective link
theory~\cite{MR91i:57001}.

\bibliographystyle{amsalpha}
\bibliography{graph-theory,knot-theory,linguistics,topology,virtual-knots,publications,preparation,immersed-curves}

\end{document}